\documentclass[a4paper,12pt]{amsart}
\usepackage{amsmath,amssymb,color,inputenc,euscript,graphicx,psfrag}
\newtheorem{theorem}{Theorem}[section]
\newtheorem{proposition}[theorem]{Proposition}
\newtheorem{lemma}[theorem]{Lemma}
\newtheorem{corollary}[theorem]{Corollary}

\newtheorem{remark}[theorem]{Remark}

\newcommand\1{1\hspace{-1mm}\text{l}}
\newcommand\co{\operatorname{co}}
\newcommand\const{\operatorname{const.}}
\renewcommand\Im{\operatorname{Im}}
\renewcommand\Re{\operatorname{Re}}
\newcommand\supp{\operatorname{supp}}

\title[Three results in Dunkl analysis]
{Three results in Dunkl analysis }
\author[B. Amri, J.--Ph. Anker, M. Sifi]
{B\'echir Amri, Jean--Philippe Anker \& Mohamed Sifi}
\address{B\'echir Amri,
Universit\'e de Tunis,
Institut Pr\'eparatoire aux Etudes d'Ing\'enieurs de Tunis (IPEIT),
D\'epartement de Math\'ematiques,
1089 Montfleury Tunis,
Tunisie}
\email{bechir.amri@ipeit.rnu.tn}
\address{Jean--Philippe Anker,
Universit\'e d'Orl\'eans \& CNRS,
Laboratoire MAPMO (UMR 6628),
F\'ed\'eration Denis Poisson (FR 2964),
B\^atiment de Math\'ematiques,
B.P.~6759,
45067 Orl\'eans cedex 2,
France}
\email{anker@univ-orleans.fr}
\address{Mohamed Sifi,
Universit\'e de Tunis El Manar,
Facult\'e des Sciences de Tunis,
D\'epartement de Math\'ematiques,
2092 Tunis,
Tunisie}
\email{mohamed.sifi@fst.rnu.tn}
\thanks{Authors partially supported
by DGRST project 04/UR/15-02 and CMCU program 07G 1501}
\date{\today}
\dedicatory{In memory of Andrzej Hulanicki (1933--2008),\\
a distinguished polish mathematician, a guide and a friend,\\
who has left many orphans in Wroc{\l}aw and around the world.\\
We miss you.}
\begin{document}

\begin{abstract}
In this article,
we establish first a geometric Paley--Wiener theorem
for the Dunkl transform in the crystallographic case.
Next we obtain an optimal bound
for the $L^p\to L^p$ norm of Dunkl translations in dimension 1.
Finally we describe more precisely the support
of the distribution associated to Dunkl translations in higher dimension.
\end{abstract}

\maketitle

\section{Introduction}
\label{section_introduction}

Dunkl theory generalizes classical Fourier analysis on \,$\mathbb{R}^N$.
It started twenty years ago with Dunkl's seminal work \cite{D1}
and was further developed by several mathematicians.
See for instance the surveys \cite{R4, G}
and the references cited therein.

In this setting,
the Paley--Wiener theorem is known to hold
for balls centered at the origin.
In \cite{J2},
a Paley--Wiener theorem was conjectured
for convex neighborhoods of the origin,
which are invariant under the underlying reflection group,
and was partially proved.
Our first result in Section \ref{section_PaleyWiener}
is a proof of this conjecture in the crystallographic case,
following the third approach in \cite{J2}.

Generalized translations were introduced in \cite{R2}
and further studied in \cite{T2, R5, TX}.
Apart from their abstract definition,
we lack precise information,
in particular about their integral representation
\begin{equation*}
(\tau_xf)(y)=\int_{\mathbb{R}^N}\hspace{-1mm}f(z)\,d\gamma_{x,y}(z)\,,
\end{equation*}
which was conjectured in \cite{R2}
and established in few cases,
for instance in dimension \,$N\!=\!1$
\,or when $f$ is radial.
Our second result in Section \ref{section_bounds}
is an optimal bound for the integral
$$
\int_{\,\mathbb{R}}|d\gamma_{x,y}(z)|
$$
in dimension \,$N\!=\!1$\,,
\,improving upon earlier results in \cite{R1, TX}.
Our bound depends on the multiplicity \,$k\!\ge\!0$
\,and tends from below to $\sqrt{2\,}$,
as \,$k\to+\infty$\,.
Our third result in Section \ref{section_support}
deals with the support of the distribution \,$\gamma_{x,y}$ \,in higher dimension,
that we determine rather precisely in the crystallographic case.

\section{Background}
\label{section_background}

In this section,
we recall some notations and results in Dunkl theory
and we refer for more details
to the articles \cite{D1, J1}
or to the surveys \cite{R4, G}.

Let $G\!\subset\!\text{O}(\mathbb{R}^N)$
be a finite reflection group associated to a reduced root system $R$
and $k:R\rightarrow[0,+\infty)$  a $G$--invariant function
(called multiplicity function).
Let $R^+$ be a positive root subsystem,
\,$\Gamma_{\!+}$ the corresponding open positive chamber,
\,$\overline{\Gamma_{\!+}}$ its closure,
\,$\overline{\Gamma^+}\hspace{-1mm}
=\!\sum_{\,\alpha\in R^+}\mathbb{R}_+\alpha$
the dual cone,
and let us denote by \,$x_+$
the intersection point of any orbit $G.x$ \,in \,$\mathbb{R}^N$
with \,$\overline{\Gamma_{\!+}}$\,.

The Dunkl operators \,$T_\xi$ on $\mathbb{R}^N$ are
the following $k$--de\-for\-ma\-tions of directional derivatives $\partial_\xi$
by difference operators\,:
$$\textstyle
T_\xi f(x)=\partial_\xi f(x)
+\sum_{\,\alpha\in R^+}\!k(\alpha)\,\langle\alpha,\xi\rangle\,
\frac{f(x)-f(\sigma_\alpha.\,x)}{\langle\alpha,\,x\rangle}\,,
$$
where \,$\sigma_\alpha.\,x=
x-\frac{\langle\alpha,\,x\rangle}{2\,|\alpha|^2}\,\alpha$
\,denotes the reflection
with respect to the hyperplane orthogonal to $\alpha$.
The Dunkl operators are antisymmetric
with respect to the measure $w(x)\,dx$
with density
$$\textstyle
w(x)=\,\prod_{\,\alpha\in R^+}|\,\langle\alpha,x\rangle\,|^{\,2\,k(\alpha)}\,.
$$

The operators $\partial_\xi$ and $T_\xi$
are intertwined by a Laplace--type operator
\begin{eqnarray}\label{vk}
V\hspace{-.25mm}f(x)\,
=\int_{\mathbb{R}^N}\hspace{-1mm}f(y)\,d\mu_x(y)
\end{eqnarray}
associated to a family of compactly supported probability measures
\,$\{\,\mu_x\,|\,x\!\in\!\mathbb{R}^N\hspace{.25mm}\}$\,.
Specifically, \,$\mu_x$ is supported in the the convex hull
$$
C^{\,x}=\,\co(G.x)\,.
$$

For every $\lambda\!\in\!\mathbb{C}^N$\!,
the simultaneous eigenfunction problem
\begin{equation*}
T_\xi f=\langle\lambda,\xi\rangle\,f
\qquad\forall\;\xi\!\in\!\mathbb{R}^N
\end{equation*}
has a unique solution $f(x)\!=\!E(\lambda,x)$
such that $E(\lambda,0)\!=\!1$,
which is given by
\begin{equation}\label{EV}
E(\lambda,x)\,
=\,V(e^{\,\langle\lambda,\,.\,\rangle})(x)\,
=\int_{\mathbb{R}^N}\hspace{-1mm}e^{\,\langle\lambda,y\rangle}\,d\mu_x(y)
\qquad\forall\;x\!\in\!\mathbb{R}^N.
\end{equation}
Furthermore \,$\lambda\mapsto E(\lambda,x)$ \,extends to
a holomorphic function on \,$\mathbb{C}^N$
and the following estimate holds\,:
\begin{equation*}
|E(\lambda,x)|\leq e^{\,\langle(\Re \lambda)_+,\,x_+\rangle}
\quad\forall\;\lambda\!\in\!\mathbb{C}^N,\;\forall\;x\!\in\!\mathbb{R}^N.
\end{equation*}

In dimension $N\!=\!1$,
these functions can be expressed in terms of Bessel functions.
Specifically,
$$\textstyle
E(\lambda,x)=j_{k-\frac12}(\lambda\,x)
+\frac{\lambda\,x}{2\hspace{.25mm}k+1}\,j_{k+\frac12}(\lambda\,x)\,,
$$
where 
$$\textstyle
j_\nu(z)\,=\;\Gamma(\nu\!+\!1)\;
{\displaystyle\sum\nolimits_{\,n=0}^{+\infty}}\;
\frac{(-1)^n}{n\hspace{.25mm}!\,\Gamma(\nu+n+1)}\;
\bigl(\frac z2)^{2n}
$$
are normalized Bessel functions.

The Dunkl transform  is defined on $L^1(\mathbb{R}^N\!,w(x)dx)$ by
$$
\mathcal{D}f(\xi)={\textstyle\frac1c}
\int_{\mathbb{R}^N}\!f(x)\,E(-i\,\xi,x)\,w(x)\,dx\,,
$$
where
$$
c\,=\int_{\mathbb{R}^N}\!e^{-\frac{|x|^2}2}\,w(x)\,dx\,.
$$
We list some known properties of this transform\,:
\begin{itemize}
\item[(i)]
The Dunkl transform is a topological automorphism
of the Schwartz space $\mathcal{S}(\mathbb{R}^N)$.
\item[(ii)]
(\textit{Plancherel Theorem\/})
The Dunkl transform extends to
an isometric automorphism of $L^2(\mathbb{R}^N\!,w(x)dx)$.
\item[(iii)]
(\textit{Inversion formula\/})
For every $f\!\in\!\mathcal{S}(\mathbb{R}^N)$,
and more generally for every $f\!\in\!L^1(\mathbb{R}^N\!,w(x)dx)$
such that $\mathcal{D}f\!\in\!L^1(\mathbb{R}^N\!,w(\xi)d\xi)$,
we have
$$
f(x)=\mathcal{D}^2\!f(-x)\qquad\forall\;x\!\in\!\mathbb{R}^N.
$$
\item[(iv)]
(\textit{Paley--Wiener theorem\/})
The Dunkl transform is a linear isomorphism between
the space of smooth functions $f$ on \,$\mathbb{R}^N$
with \,$\supp f\!\subset\hspace{-.5mm}\overline{B(0,R)}$
\,and the space of entire functions $h$ on \,$\mathbb{C}^N$
such that
\begin{equation}\label{PaleyWiener0}\textstyle
\sup_{\,\xi\in\mathbb{C}^N}\,(1\!+\!|\xi|)^M\,e^{-R\,|\Im\xi\,|}\,|h(\xi)|<+\infty
\qquad\forall\;M\!\in\!\mathbb{N}\,.
\end{equation}
\end{itemize}

\section{A geometric Paley--Wiener theorem}
\label{section_PaleyWiener}

In this section,
we prove a geometric version of the Paley--Wiener theorem,
which was looked for in \cite{J2, T2, J3},
under the assumption that $G$ is crystallographic.
The proof consists merely in resuming the third approach in \cite{J2}
and applying it to the convex sets considered in \cite{A1, A2, A3, A4}
instead of the convex sets considered in \cite{O1}.
Recall that the second family consists of the convex hulls
$$
C^\Lambda=\co(G.\Lambda)
$$
of $G$--orbits $G.\Lambda$ in $\mathbb{R}^N$,
while the first family consists of the polar sets
$$
C_\Lambda=\{\,x\!\in\!\mathbb{R}^N\mid
\langle x,g.\Lambda\rangle\!\le\!1\hspace{2mm}\forall\;g\!\in\!G\,\}\,.
$$

\begin{figure}[ht]
\begin{center}
\psfrag{Lambda}[c]{\color{red}$\Lambda$}
\includegraphics[width=55mm]{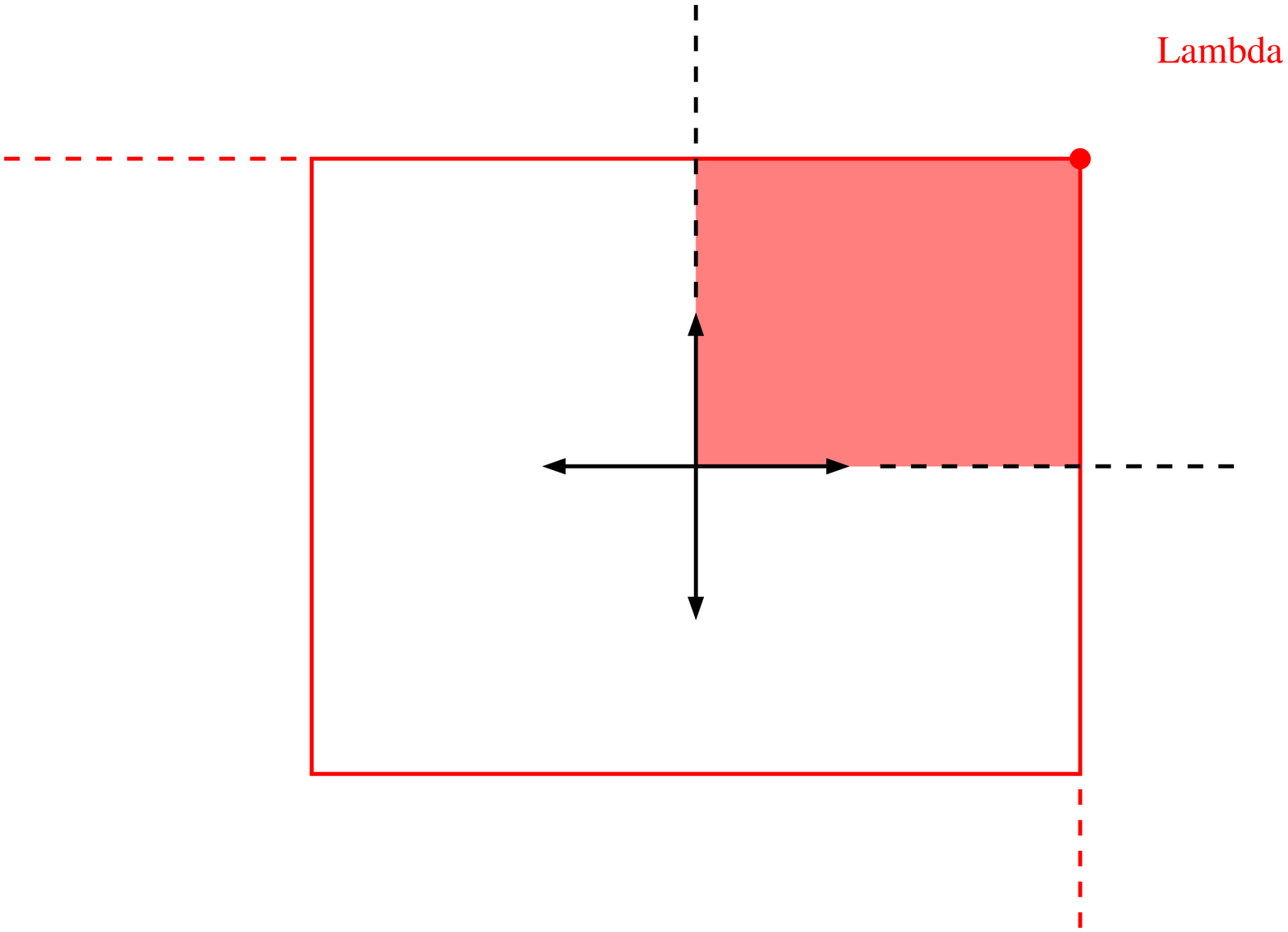}
\hspace{15mm}
\includegraphics[width=40mm]{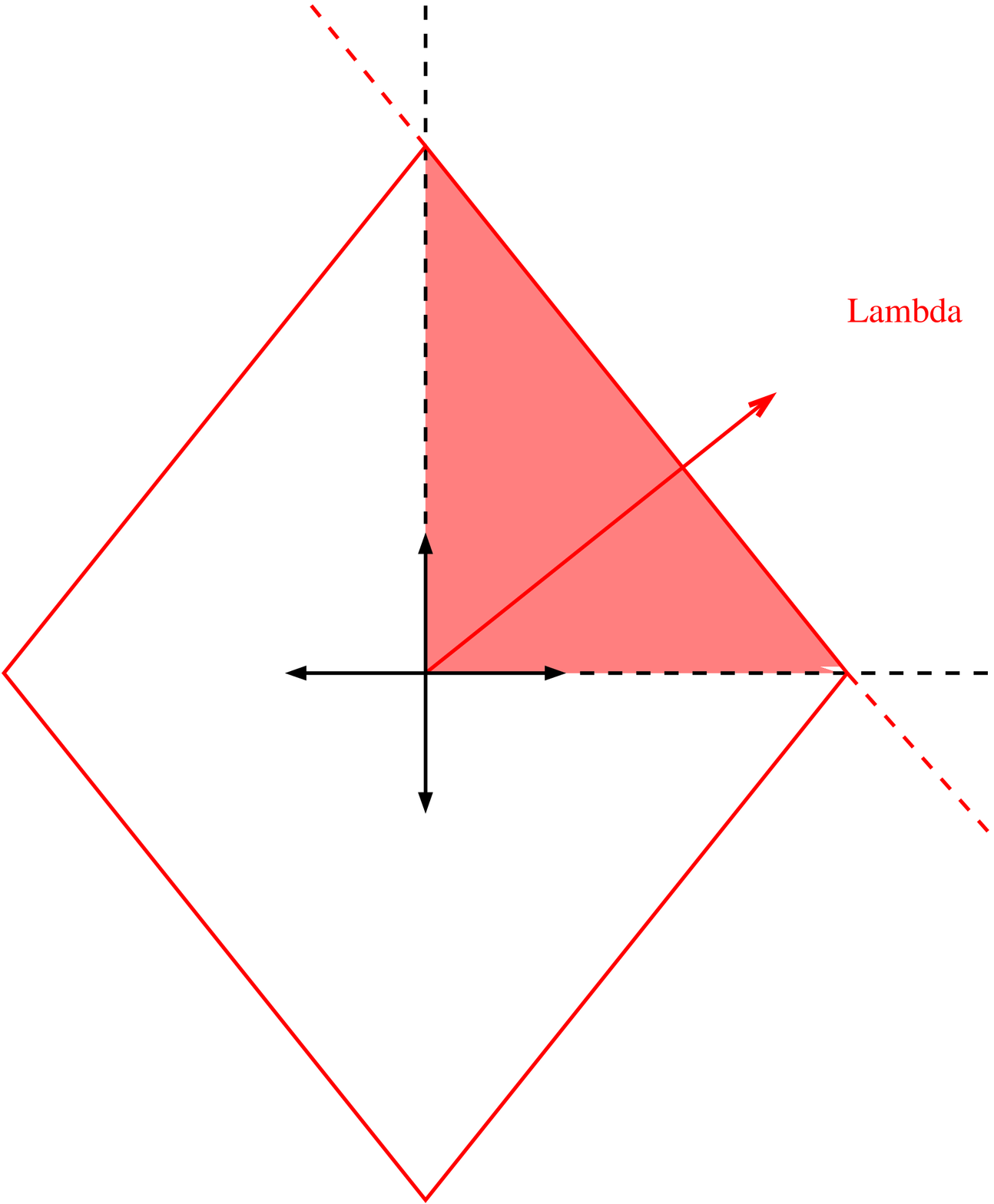}
\caption{The sets \,$C^\Lambda$ \,and \,$C_\Lambda$
\,for the root system \,$A_1\!\times\!A_1$}
\end{center}
\end{figure}

\begin{figure}[ht]
\begin{center}
\psfrag{Lambda}[c]{\color{red}$\Lambda$}
\includegraphics[width=50mm]{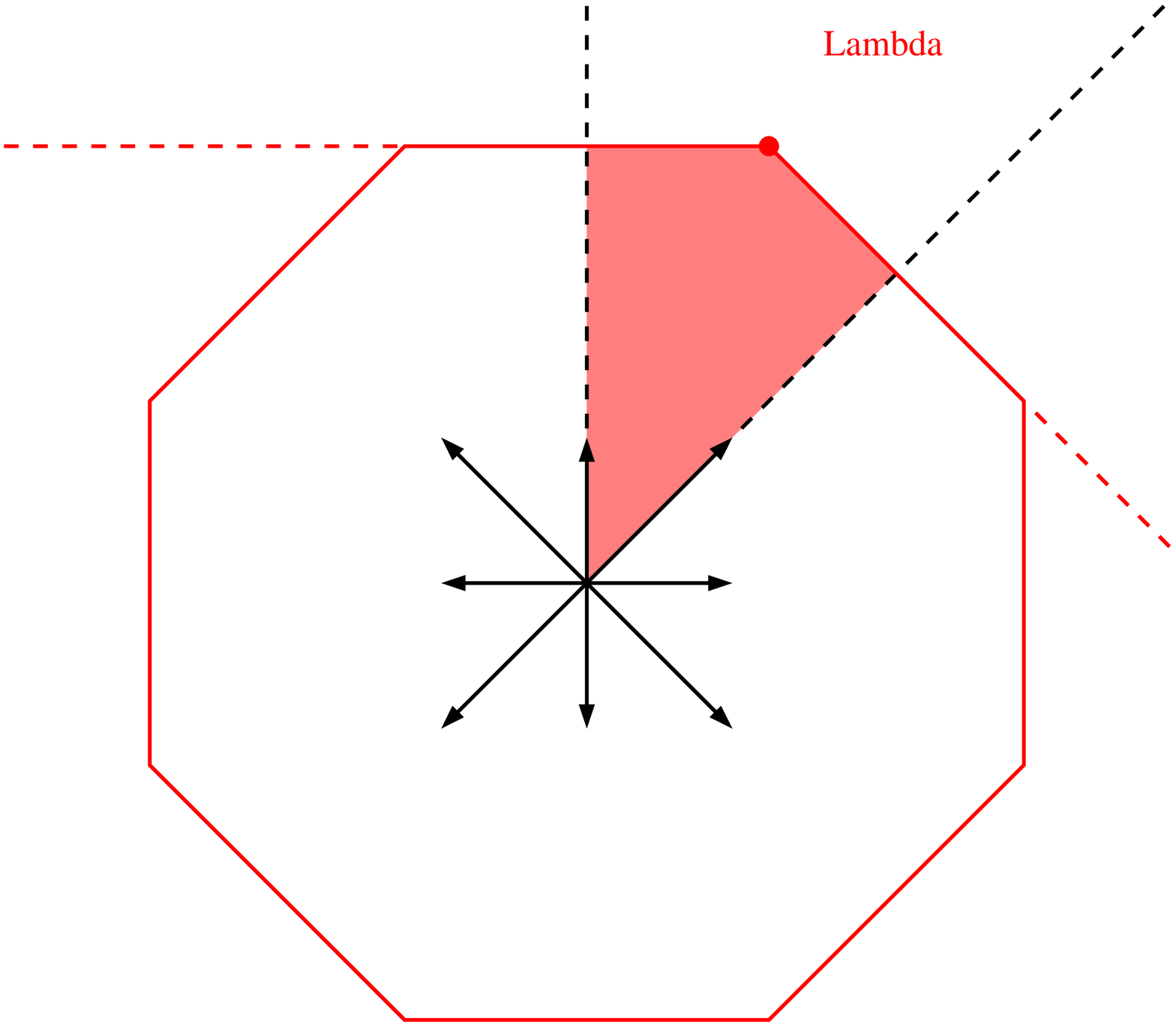}
\hspace{10mm}
\includegraphics[width=50mm]{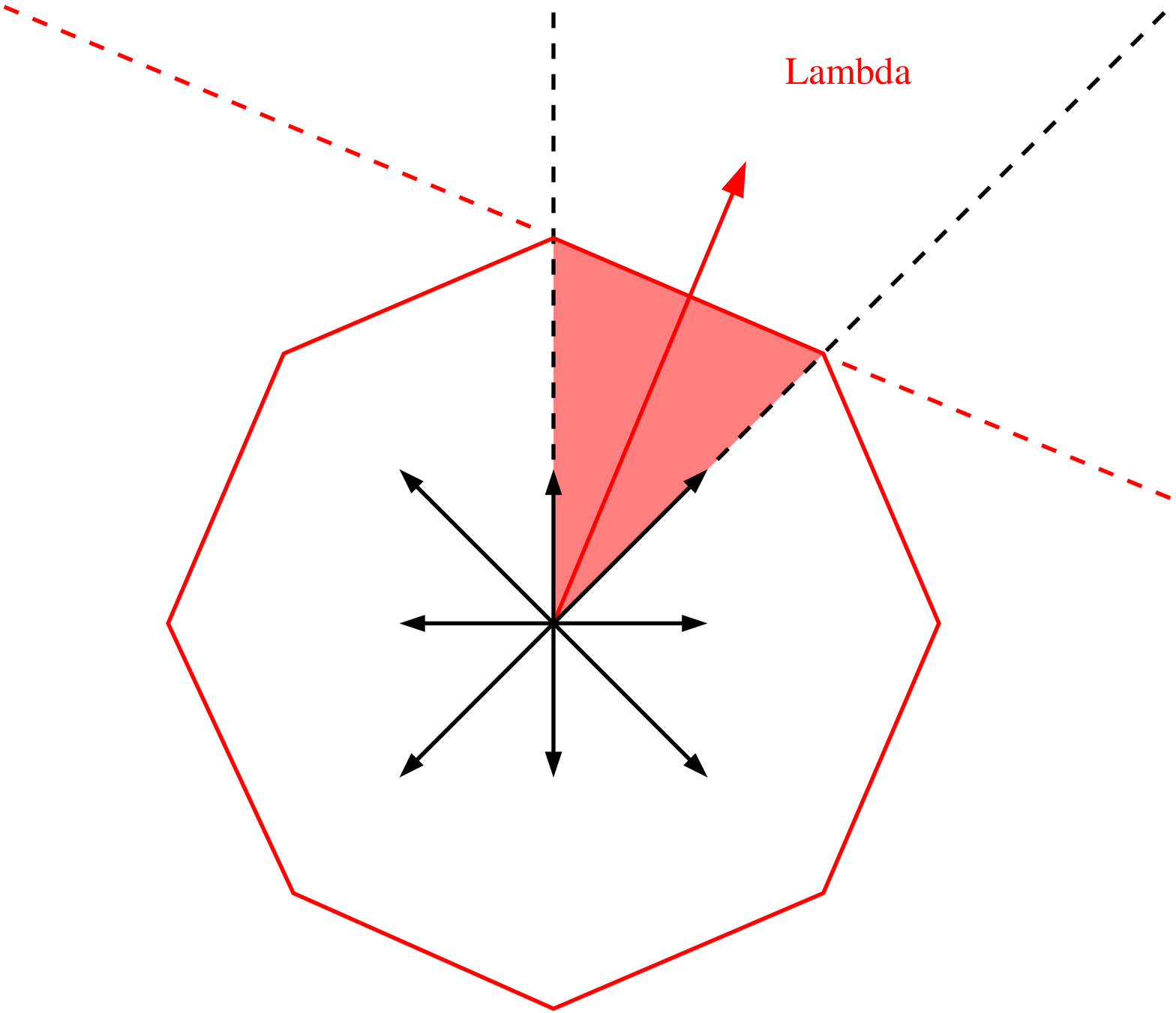}
\caption{The sets \,$C^\Lambda$ \,and \,$C_\Lambda$
\,for the root system \,$B_2$}
\end{center}
\end{figure}

Before stating the geometric Paley--Wiener theorem,
let us make some remarks about the sets $C^\Lambda$ and $C_\Lambda$.
Firstly, they are convex, closed, $G$--invariant
and the following inclusion holds\,:

\centerline{$
C^\Lambda\,\subset\,|\Lambda|^2\,C_\Lambda\,.
$}\vspace{1mm}

\noindent
Secondly, we may always assume that
\,$\Lambda\!=\!\Lambda_+$ belongs to
the closed positive chamber \,$\overline{\Gamma_{\!+}}$
and, in this case, we have
\begin{align*}
C^\Lambda\cap\,\overline{\Gamma_{\!+}}\,
&=\;\textstyle\overline{\Gamma_{\!+}}\,\cap
\bigl(\,\Lambda-\overline{\Gamma^+}\bigr)\,,\\
C_\Lambda\cap\,\overline{\Gamma_{\!+}}\,
&=\,\{\,x\!\in\!\overline{\Gamma_{\!+}}\mid
\langle\Lambda,x\rangle\le1\,\}\,.
\end{align*}
Thirdly, on one hand,
every $G$--invariant convex subset in $\mathbb{R}^N$
is a union of sets $C^\Lambda$
while, on the other hand,
every $G$--invariant closed convex subset in $\mathbb{R}^N$
is an intersection of sets $C_\Lambda$.
For instance,
$$
\overline{B(0,R)}\,
=\;\bigcup\nolimits_{\,|\Lambda|=R}C^{\Lambda}\,
=\;\bigcap\nolimits_{\,|\Lambda|=R^{-1}}\hspace{-.5mm}C_{\Lambda}\,.
$$
Fourthly, we shall say that
$\Lambda\!\in\!\overline{\Gamma_+}$ is {\it admissible\/}
if the following equivalent conditions are satisfied\,:

\begin{itemize}

\item[(i)]
\,$\Lambda$ \hspace{.5mm}has nonzero projections
in each irreducible component of
\hspace{.25mm}$(\mathbb{R}^N\hspace{-.5mm},
\hspace{-.25mm}R)$\hspace{.25mm},

\item[(ii)]
\,$C^\Lambda$ is a neighborhood of the origin,

\item[(iii)]
\,$C_\Lambda$ is bounded.

\end{itemize}

\noindent
In this case, we may consider the gauge
$$\textstyle
\chi_\Lambda(\xi)
=\max_{\,x\in C_\Lambda}\langle x,\xi\rangle
=\min\,\{\,r\!\in[\hspace{.5mm}0,+\infty)\mid\xi\!\in\!r\,C^\Lambda\,\}
$$
 on \hspace{.5mm}$\mathbb{R}^N$.
 
\begin{theorem}\label{PaleyWiener1}
Assume that \,$\Lambda\!\in\!\overline{\Gamma_+}$ is admissible.
Then the Dunkl transform is a linear isomorphism
between the space of smooth functions \hspace{.5mm}$f$ on \,$\mathbb{R}^N$
with \,$\supp f\!\subset\hspace{-.5mm}C_\Lambda$
and the space of entire functions \hspace{.5mm}$h$ on \,$\mathbb{C}^N$
such that
\begin{equation}\label{PW1}\textstyle
\sup_{\,\xi\in\mathbb{C}^N}
(1\!+\!|\xi|)^M\,e^{-\chi_\Lambda(\Im\xi)}\,|h(\xi)|
<+\infty
\qquad\forall\;M\!\in\!\mathbb{N}\,.
\end{equation}
\end{theorem}

\begin{proof}
Following \cite{J2},
this theorem is first proved in the trigonometric case,
which explains the restriction to crystallographic groups,
and next obtained in the rational case by passing to the limit.
The proof of Theorem \ref{PaleyWiener1} in the trigonometric case
is similar to the proof of the Paley--Wiener Theorem in \cite{O1, O2},
and actually to the initial proof of Helgason
for the spherical Fourier transform
on symmetric spaces of the noncompact type.
This was already observed in \cite{S}
and will be developed below for the reader's convenience.
The limiting procedure, as far as it is concerned,
is described thoroughly in \cite{J2}
and needs no further explanation.

Thus assume that
$h$ is an entire function on $\mathbb{C}^N$ satisfying \eqref{PW1}
and, by resuming the proof of \cite[Theorem 8.6\,(2)]{O1},
let us show that its inverse Cherednik transform
\begin{equation}\label{ICT}
f(x)\hspace{.25mm}=\hspace{.5mm}\const\int_{\mathbb{R}^N}\!
h(\xi)\,\widetilde E(i\hspace{.25mm}\xi,x)\,\widetilde w(\xi)\,d\xi
\end{equation}
vanishes outside $C_\Lambda$.
Firstly, one may restrict by $G$--equivariance to $x\!=\!g_0.x_+$,
where \hspace{.25mm}$x_+\hspace{-1mm}\in\!
\Gamma_+\hspace{-1mm}\smallsetminus \!C_\Lambda$
and $g_0$ denotes the longest element in $G$,
which interchanges $\Gamma_+$ and $-\Gamma_+$.
Secondly, by expanding
\begin{equation*}
\Bigl\{\,\prod\nolimits_{\,\alpha\in R^+}
(\hspace{.25mm}\langle\check\alpha,\xi\rangle\!
-\!k_\alpha\hspace{.25mm})\hspace{.25mm}\Bigr\}\;
\widetilde E(\xi,x)\,
=\,\sum\nolimits_{\,g\in G}\sum\nolimits_{\,q\in Q^+}\!
{\mathbf c}(-g.\xi)\,\widetilde E_q(g,g.\xi)\,
e^{\hspace{.25mm}\langle\hspace{.25mm}g.\xi
+\varrho\hspace{.2mm}+\hspace{.1mm}q,
\hspace{.5mm}x\hspace{.25mm}\rangle}
\end{equation*}
\eqref{ICT} becomes
\begin{equation*}
f(x)\hspace{.25mm}=\hspace{.5mm}\const\hspace{.5mm}
\sum\nolimits_{\,g\in G}\det g\,\sum\nolimits_{\,q\in Q^+}f_{g,q}(x)\,
e^{\hspace{.25mm}\langle\hspace{.25mm}\varrho\hspace{.2mm}+\hspace{.1mm}q,
\hspace{.5mm}x\hspace{.25mm}\rangle}\,,
\end{equation*}
where
\begin{equation}\label{coefficients}
f_{g,q}(x)\,=\int_{\mathbb{R}^N}\!h(g^{-1}\hspace{-.5mm}.\hspace{.25mm}\xi)\,
\widetilde E_q(g,i\hspace{.25mm}\xi)\,
e^{\hspace{.5mm}i\hspace{.25mm}\langle\hspace{.25mm}\xi,
\hspace{.5mm}x\hspace{.25mm}\rangle}\,
\Bigl\{\hspace{.5mm}\prod\nolimits_{\,\alpha\in R^+}\hspace{-1mm}\textstyle
\frac{\Gamma(\hspace{.25mm}i\hspace{.25mm}\langle\hspace{.25mm}\check\alpha,
\hspace{.5mm}\xi\hspace{.25mm}\rangle\hspace{.25mm}
+\hspace{.25mm}k_\alpha\hspace{.25mm})}
{\Gamma(\hspace{.25mm}i\hspace{.25mm}\langle\hspace{.25mm}\check\alpha,
\hspace{.5mm}\xi\hspace{.25mm}\rangle\hspace{.25mm}
+\hspace{.25mm}1\hspace{.25mm})}\hspace{.25mm}\Bigr\}\,d\xi\,.
\end{equation}
Thirdly, one shows that all expressions \eqref{coefficients} vanish,
by shifting the contour of integration from \,$\mathbb{R}^N$ to
\,$\mathbb{R}^N\!+i\hspace{.5mm}t\hspace{.5mm}
g_0\hspace{.1mm}.\hspace{.2mm}\Lambda$
\hspace{.25mm}with \,$t\!>\!0$\hspace{.25mm},
which produces an exponential factor \,$e^{-c\hspace{.25mm}t}$
\,with \,$c\hspace{-.25mm}=\!\langle\Lambda,x_+\rangle\!-\!1\hspace{-.5mm}
>\hspace{-.5mm}0$\hspace{.5mm},
and by letting \,$t\!\to\!+\infty$\,.
\end{proof}

Since every $G$--invariant convex compact neighborhood
of the origin in \hspace{.5mm}$\mathbb{R}^N$
is the intersection of admissible sets \hspace{.25mm}$C_\Lambda$\hspace{.25mm},
Theorem \ref{PaleyWiener1} generalizes as follows.

\begin{corollary}[Geometric Paley--Wiener Theorem]\label{PaleyWiener2}
Let $C$ be a $G$--invariant convex compact neighborhood
of the origin in $\mathbb{R}^N$
and $\chi(\xi)\!=\!\max_{\,x\in C}\langle x,\xi\rangle$ the dual gauge.
Then the Dunkl transform is a linear isomorphism
between the space \,$\mathcal{C}_C^\infty(\mathbb{R}^N)$
of smooth functions $f$ on \,$\mathbb{R}^N$
with \,$\supp f\!\subset\hspace{-.5mm}C$
\,and the space \,$\mathcal{H}_\chi(\mathbb{C}^N)$
of entire functions $h$ on \,$\mathbb{C}^N$ such that
$$\textstyle
\sup_{\,\xi\in\mathbb{C}^N}
(1\!+\!|\xi|)^M\,e^{-\chi(\Im\xi)}\,|h(\xi)|
<+\infty
\qquad\forall\;M\!\in\!\mathbb{N}\,.
$$
\end{corollary}

\begin{remark}
Notice that the Dunkl transform \,$\mathcal{D}$
always maps \,$\mathcal{C}_C^\infty(\mathbb{R}^N)$
into \,$\mathcal{H}_\chi(\mathbb{C}^N)$
and that the assumption that $G$ is crystallographic
is only used to prove that \,$\mathcal{D}$ is onto.
\end{remark}

\section{$L^p$ bounds for generalized translations in dimension $1$}
\label{section_bounds}

Dunkl translations are defined on $\mathcal{S}(\mathbb{R}^N)$ by
$$
(\tau_xf)(y)=
{\textstyle\frac1c}\int_{\mathbb{R}^N}\hspace{-1mm}
\mathcal{D}f(\xi)\,E(i\hspace{.25mm}\xi,x)\,E(i\hspace{.25mm}\xi,y)\,w(\xi)\,d\xi
\qquad\forall\;x,y\!\in\!\mathbb{R}^N.
$$
They have an explicit integral representation \cite{R1} in dimension $N\!=\!1$\,:
$$
(\tau_xf)(y)=\int_{\mathbb{R}}f(z)\,d\gamma_{x,y}(z)\,,
$$
where
\begin{equation}\label{gm}
d\gamma_{x,y}(z)\,=\,\begin{cases}
\,\gamma(x,y,z)\,|z|^{2k}\,dz
&\text{if \,}x,y\!\in\!\mathbb{R}^*\\
\,d\delta_y(z)
&\text{if \,}x\!=\!0\\
\,d\delta_x(z)
&\text{if \,}y\!=\!0
\end{cases}
\end{equation}
is a signed measure such that
\,$\displaystyle\int_{\mathbb{R}}d\gamma_{x,y}(z)=1$\,.
Specifically,
$$\textstyle
\gamma(x,y,z)\,=\,d\;\sigma(x,y,z)\;\rho(|x|,|y|,|z|)\;\1_{\,I_{|x|,|y|}}(|z|)
\qquad\forall\;x,y,z\in\mathbb{R}^*,
$$
where
\begin{equation*}
\begin{aligned}
d&=\textstyle
\frac{\Gamma(k+\frac12)}{\sqrt{\pi}\,\Gamma(k)\vphantom{\frac12}}\,,\\
\sigma(x,y,z)
&=\textstyle
1-\frac{x^2+\,y^2\,-z^2}{2\,x\,y}
+\frac{z^2+\,y^2\,-x^2}{2\,z\,y}
+\frac{x^2+\,z^2\,-y^2}{2\,x\,z}\\
&=\textstyle
\frac{(z+x+y)\,(z+x-y)\,(z-x+y)}{2\,x\,y\,z}
\qquad\forall\;x,y,z\in\mathbb{R}^*,\\
\rho(a,b,c)&=\textstyle
\frac{\{\,c^2-\,(a-b)^2\}^{k-1}\,
\{\,(a+b)^2-\,c^2\}^{k-1}\vphantom{\frac12}}
{(\,2\,a\,b\,c\,)^{\,2k-1}\vphantom{\frac12}}\\
&=\textstyle
\frac{(\,2\,b^2c^2+\,2\,a^2c^2+\,2\,a^2b^2
-\,a^4-\,b^4-\,c^4\vphantom{\frac12})^{\,k-1}}
{(\,2\,a\,b\,c\,)^{\,2k-1}\vphantom{\frac12}}
\qquad\forall\;a,b,c>0\,,
\end{aligned}
\end{equation*}
and \,$I_{a,b}$ \,denotes the interval \,$[\,|\,a\!-\!b\,|,a\!+\!b\,]$\,.
Notice the symmetries
\begin{equation}\label{gamma}
\gamma(x,y,z)\,=\,\begin{cases}
\,\gamma(y,x,z)\,,\\
\,\gamma(-x,-y,-z)\,,\\
\,\gamma(-z,y,-x)=\gamma(x,-z,-y)\,.
\end{cases}
\end{equation}

\begin{proposition}\label{measure}
The following inequality holds, for every \,$x,y\!\in\!\mathbb{R}$\,:
\begin{equation}\label{A}
\int_{\,\mathbb{R}}\,\bigl|d\gamma_{x,y}(z)\bigr|\leq A_k=
\textstyle\sqrt{2}\,
\frac{\{\,\Gamma(k+\frac12)\}^2}{\Gamma(k+\frac14)\,\Gamma(k+\frac34)}\,.
\end{equation}
Actually there is equality if \,$x\!=\!y\!\in\!\mathbb R^*$.
Moreover \,$A_k\!\overset<\longrightarrow\!\sqrt{2}$
\,as \,$k\!\to\!+\infty$\,.
\end{proposition}

\begin{remark}
This result improves earlier bounds obtained in \cite{R1} and \cite{TX},
which were respectively \,$4$ and \,$3$\,.

\end{remark}

\begin{proof}
Let $x,y\!\in\!\mathbb{R}^*$.
\vskip2mm

\noindent
\textit{Case 1}\,: Assume that \,$x\,y\!<\!0$\,.
Then \,$|\,|x|\!-\!|y|\,|=|\,x\!+\!y\,|$ \,and \,$|x|\!+\!|y|=|\,x\!-\!y\,|$\,,
hence \,$\sigma(x,y,z)\,\1_{\,I_{|x|,|y|}}(|z|)
=\frac{z\,+\,x\,+\,y}z\,\frac{(x-y)^2-\,z^2}{-\,2\,x\,y}\,\1_{\,I_{|x|,|y|}}(|z|)$
\,and \,$\gamma_{x,y}$ \,are positive.
Thus
$$
\int_{\mathbb{R}}|d\gamma_{x,y} (z)|
=\int_{\mathbb{R}}d\gamma_{x,y} (z)
=1\,.
$$
\textit{Case 2}\,: Assume that \,$x\,y\!>\!0$\,.
By symmetry, we may reduce to \,$0\!<\!x\!\leq\!y$\,.
Then
\begin{equation*}
\begin{aligned}
\int_{\,\mathbb{R}}|d\gamma_{x,y} (z)|
&\,=\int_{-\infty}^{\,0}\hspace{-1mm}|d\gamma_{x,y}(z)|\,
+\int_{\,0}^{+\infty}\hspace{-1.5mm}|d\gamma_{x,y} (z)|\\
&\,=\,2\,d\int_{y-x}^{\,y+x}\hspace{-1.5mm}
\textstyle
\frac{x\,+\,y}{2\,x\,y\,z}
\bigl(\frac{z^2-\,x^2-\,y^2\,+\,2\,x\,y}{2\,x\,y\,z}\bigr)^k\,
\bigl(\frac{x^2+\,y^2\,+\,2\,x\,y\,-\,z^2}{2\,x\,y\,z}\bigr)^{k-1}\,
z^{2k}\,dz\,.
\end{aligned}
\end{equation*}
After performing the change of variables
\,$z=\!\sqrt{\,x^2\!+y^2\!-2\,x\,y\cos\theta\,}$
\,and setting \,$y\!=\!s\,x$\,,
we get
\begin{equation}\label{gam}
\int_{\,\mathbb{R}}|d\gamma_{x,y}(z)|\,
=\,{\textstyle\frac{\Gamma(k+\frac12)}{\sqrt{\pi}\,\Gamma(k)}}\,
(1\!+\!s)\int_{\,0}^{\,\pi}\!\textstyle
\frac{(1\,-\,\cos\theta)\,\sin^{2k-1}\!\theta}
{\sqrt{1\,+\,s^2-\,2\,s\,\cos\theta}}\,d\theta\,.
\end{equation}
Denote by $F(s)$ the right hand side of \eqref{gam}.
Since
\begin{equation*}
F'(s)=\,{\textstyle\frac{\Gamma(k+\frac12)}{\sqrt{\pi}\,\Gamma(k)}}\,
(1\!-\!s)\int_{\,0}^{\,\pi}\hspace{-1mm}\textstyle
\frac{\sin^{2k+1}\!\theta\vphantom{\big|}}
{(1\,+\,s^2-\,2\,s\,\cos\theta)^{\frac32}}\,d\theta
\end{equation*}
is nonpositive,
\,$F(s)$ \,is a decreasing function on $[\,1,+\infty\,)$,
which reaches its maximum at \,$s\!=\!1$\,.
Let us compute it\,:
\begin{equation*}
\begin{aligned}
A_k=\,F(1)
&=\,{\textstyle\frac{\sqrt{2}\,\Gamma(k+\frac12)}{\sqrt{\pi}\,\Gamma(k)}}
\int_{\,0}^{\,\pi}\!(1\!-\!\cos\theta)^{k-\frac12}\,
(1\!+\!\cos\theta)^{k-1}\,\sin\theta\;d\theta\\
&=\,2^{2k}\,{\textstyle\frac{\Gamma(k+\frac12)}{\sqrt{\pi}\,\Gamma(k)}}
\int_{\,0}^{\,1}t^{k-\frac12}\,(1\!-\!t)^{k-1}\,dt\\
&=\textstyle
\,2^{2k}\,{\textstyle\frac{\Gamma(k+\frac12)}{\sqrt{\pi}\,\Gamma(k)}}\;
B(k\!+\!\frac12,k)=\textstyle
\,2^{2k}\,
\frac{\{\,\Gamma(k+\frac12)\}^2}
{\sqrt{\pi}\;\Gamma(2k+\frac12)}\,
=\,\sqrt{2}\;
\frac{\Gamma(k+\frac12)}{\Gamma(k+\frac14)}\,
\frac{\Gamma(k+\frac12)}{\Gamma(k+\frac34)}\;,
\end{aligned}
\end{equation*}
after performing the change of variables \,$t\!=\!\frac{1-\cos\theta}2$
\,and using standard properties of the beta and gamma functions.
\vspace{2mm}

Finally let us show that
\,$A_k\!\overset<\longrightarrow\!\sqrt{2}$
\,as \,$k\!\to\!+\infty$\,.
Write
$$\textstyle
A_k=\sqrt{2}\;\frac{G(k+\frac14)}{G(k+\frac12)}\;,
\quad\text{where}\quad
G(u)=\frac{\Gamma(u+\frac14)}{\Gamma(u)}
\quad\forall\;u\!>\!0\,.
$$
Since the logarithmic derivative $\frac{\Gamma'}\Gamma$
of the gamma function
is a strictly increasing analytic function on $(0,+\infty)$,
the logarithmic derivative
$$\textstyle
\frac{G'(u)}{G(u)}
=\frac{\Gamma'(u+\frac14)}{\Gamma(u+\frac14)}
-\frac{\Gamma'(u)}{\Gamma(u)}
$$
is positive.
Hence $G$ is an strictly increasing function and $A_k\!<\!\sqrt{2}$\,.
On the other hand, using Stirling's formula
$$\textstyle
\Gamma(u)\,\sim\,\sqrt{2\pi}\;u^{u-\frac12}\;e^{-u}
\quad\text{as}\quad u\to+\infty\,,
$$
we get \,$G(k\!+\!\frac14)\sim G(k\!+\!\frac12)$
\,hence \,$A_k\to\sqrt{2}$\,,
\,as \,$k\to+\infty$\,.
\end{proof}

\noindent
As a first consequence,
we obtain the $L^1\!\to\!L^1$ operator norm
of Dunkl translations in dimension $N\!=\!1$\,.

\begin{corollary}\label{L1}
Let \,$x\!\in \mathbb{R}^*$.
Then \,$\tau_x$ is a bounded operator
on \,$L^1(\mathbb{R},|x|^{2k}dx)$,
with \,$\|\tau_x\hspace{.25mm}\|_{L^1\to L^1}=A_k$\,.
\end{corollary}

\begin{proof}
The inequality 
\,$\|\tau_x\hspace{.25mm}\|_{L^1\to L^1}\!\le\hspace{-.25mm}A_k$
\,follows from \eqref{A}, together with \eqref{gamma},
and it remains for us to prove the converse inequality.
By symmetry, we may assume that $x\!>\!0$\,.
Since
$$
A_k=\,\lim\nolimits_{\,y\to x}\int_{\,\mathbb{R}}|\gamma(x,y,z)|\,|z|^{2k}\,dz\,,
$$
for every \,$0\!<\!\varepsilon\!<\!A_k$\hspace{.25mm},
there exists \,$0\!<\!\eta\!<\!x$ \,such that,
for every \,$y\!\in\hspace{-.25mm}
[\hspace{.25mm}x\!-\!\eta,x\!+\!±\eta\hspace{.25mm}]$\hspace{.25mm},
\begin{equation}\label{crucial}
\int_{\,\mathbb{R}}|\gamma(x,y,z)|\,|z|^{2k}\,dz\hspace{.5mm}
>\hspace{.25mm}A_k\hspace{-.5mm}-\hspace{-.25mm}\varepsilon\,.
\end{equation}
Let $f$ be a nonnegative measurable function on \,$\mathbb{R}$
\,such that
$$
\supp f\!\subset\hspace{-.25mm}
[-\hspace{.25mm}x\!-\!\eta,-\hspace{.25mm}x\!+\!±\eta\hspace{.5mm}]
\quad\text{and}\quad
\bigl\|\hspace{.25mm}f\hspace{.25mm}\bigr\|_{L^1}\hspace{-.5mm}
=\hspace{-.5mm}\displaystyle\int_{\,\mathbb{R}}\!f(z)\,|z|^{2k}dz=1\,.
$$
Since
$$
\begin{cases}
\,\gamma(x,y,z)\ge0
&\forall\;y\!<\!0\,,\;\forall\;z\!<\!0\,,\\
\,\gamma(x,y,z)\le0
&\forall\;y\!>\!0\,,\;\forall\;z\!<\!0\,,
\end{cases}
$$
we have
$$
\bigl|(\tau_xf)(y)\bigr|=
\int_{-x-\eta}^{-x+\eta}\hspace{-1mm}
f(z)\,|\gamma(x,y,z)|\,|z|^{2k}\,dz\,.
$$
Hence, using \eqref{gamma} and \eqref{crucial},
$$
\bigl\|\hspace{.5mm}\tau_xf\hspace{.5mm}\bigr\|_{L^1}
=\int_{\,\mathbb{R}}\hspace{.5mm}\bigl|(\tau_xf)(y)\bigr|\,|y|^{2k}\,dy\,
=\int_{-x-\eta}^{-x+\eta}\Bigl\{\,
\int_{\,\mathbb{R}}|\gamma(x,-z,-y)|\,|y|^{2k}\,dy\,
\Bigr\}\,f(z)\,|z|^{2k}\,dz
$$
is bounded from below by
\,$A_k\hspace{-.5mm}-\hspace{-.25mm}\varepsilon$\,.
Consequently
\,$\|\tau_x\hspace{.25mm}\|_{L^1\to L^1}\!
\ge A_k\hspace{-.5mm}-\hspace{-.25mm}\varepsilon$
\,and we conclude by letting \,$\varepsilon\!\to\!0$\,.
\end{proof}

Let us next compute the $L^2\!\to\!L^2$ operator norm of Dunkl translations.

\begin{lemma}\label{L2}
Let \,$x\!\in \mathbb{R}$\,.
Then \,$\tau_x$ is a bounded operator
on \,$L^2(\mathbb{R},|x|^{2k}dx)$,
with \,$\|\tau_x\hspace{.25mm}\|_{L^2\to L^2}=1$\,.
\end{lemma}

\begin{proof}
The proof is straightforward, via the Plancherel formula,
and generalizes to higher dimensions.
On one hand, the inequality \,$\|\tau_x\|_{L^2\to L^2}\!\le\!1$
\,follows from the estimate \,$|E(i\hspace{.25mm}\xi,x)|\!\le\!1$\,.
On the other hand, let
$$
f_\varepsilon(x)=
\varepsilon^{\hspace{.25mm}k+\frac12}\,
f(\varepsilon\hspace{.25mm}x)
$$
be a rescaled normalized function in $L^2(\mathbb{R},|x|^{2k}\hspace{.25mm}dx)$.
Then
\begin{equation*}
\|\hspace{.5mm}f_\varepsilon\hspace{.5mm}\|_{L^2}
=\|\hspace{.5mm}f\hspace{.5mm}\|_{L^2}
=1
\end{equation*}
while
\begin{equation*}
\begin{aligned}
\|\hspace{.5mm}\tau_x\hspace{.5mm}f_\varepsilon\hspace{.5mm}\|_{L^2}^{\,2}
&=\int_{\,\mathbb{R}}|E(i\hspace{.25mm}\xi,x)|^2\,
\varepsilon^{-2k-1}\,|\mathcal{D}\hspace{-.25mm}f(\varepsilon^{-1}\xi)|^2\,
|\xi|^{2k}\,d\xi\\
&=\int_{\,\mathbb{R}}|E(i\hspace{.25mm}\varepsilon\hspace{.25mm}\xi,x)|^2\,
|\mathcal{D}\hspace{-.25mm}f(\xi)|^2\,|\xi|^{2k}\,d\xi
\end{aligned}
\end{equation*}
tends to
$$
\int_{\,\mathbb{R}}
|\hspace{.25mm}\mathcal{D}\hspace{-.25mm}f(\xi)|^2\,
|\xi|^{2k}\,d\xi\hspace{.5mm}
=\|f\|_{L^2}^{\,2}
=1
$$
as \,$\varepsilon\to0$\,.
This concludes the proof of the lemma.
\end{proof}

Eventually, Corollary \ref{L1} and Lemma \ref{L2}
imply the following result, by interpolation and duality.

\begin{corollary}\label{Lp}
Let \,$x\!\in \mathbb{R}$ \,and \,$1\!\le\!p\!\le\!\infty$\,.
Then \,$\tau_x$ is a bounded operator on \,$L^p(\mathbb{R},|x|^{2k}dx)$,
with \,$\|\tau_x\hspace{.25mm}\|_{L^p\to L^p}\le
A_{\hspace{.25mm}k}^{2\,|1/p\hspace{.25mm}-1/2\hspace{.25mm}|}$\,.
\end{corollary}

\begin{remark}
In the product case,
where \,$G\!=\!\mathbb{Z}_2^N$ acts on \,$\mathbb{R}^N$,
we have
$$
\bigl\|\,\tau_x\,\bigr\|_{L^p\to L^p}\le\,
A_{\hspace{.25mm}k}^{2\,|\frac1p-\frac12|N}
$$
for every \,$x\!\in\!\mathbb{R}^N$ and \,$1\!\le\!p\!\le\!\infty$\,.
\end{remark}

\section{A support theorem for generalized translations}
\label{section_support}

As mentioned in the introduction,
we lack information about Dunkl translations in general.
In this section, we locate more precisely
the support of the distribution
$$
\langle\,\gamma_{x,y},f\,\rangle=(\tau_xf)(y)\,
$$
which is known \cite{T2} to be contained
in the closed ball of radius \,$|x|\!+\!|y|$\,.

\begin{theorem}\label{support}
{\rm (i)}
The distribution \,$\gamma_{x,y}$ \,is supported in the spherical shell
$$
\bigl\{\,z\!\in\!\mathbb{R}^N\bigm|
\bigl||x|\!-\!|y|\bigr|\le|z|\le|x|\!+\!|y|\,\bigr\}\,.
$$
{\rm (ii)}
If \,$G$ is crystallographic,
then the support of \,$\gamma_{x,y}$ \,is more precisely contained in
$$
\bigl\{\,z\!\in\!\mathbb{R}^N\bigm|
z_+\hspace{-1mm}\preccurlyeq\hspace{-.5mm}
x_+\hspace{-1mm}+\hspace{-.5mm}y_+,\,
z_+\hspace{-.25mm}\succcurlyeq
y_+\hspace{-1mm}+\hspace{-.4mm}g_0.x_+
\hspace{1mm}\text{and }\hspace{1mm}
x_+\hspace{-1mm}+\hspace{-.4mm}g_0.y_+
\,\bigr\}\,.
$$
Here \,$g_0$ denotes the longest element in $G$,
which interchanges the chambers \,$\Gamma_{\!+}$ and $-\Gamma_{\!+}$,
and \,$\preccurlyeq$ the partial order on \,$\mathbb{R}^N$
associated to the cone \,$\overline{\Gamma^+}$\,:
$$
a\preccurlyeq b
\quad\Longleftrightarrow\quad
b-a\in\overline{\Gamma^+}\,.
$$
\end{theorem}

\begin{figure}[ht]
\begin{center}
\psfrag{x+y}[r]{\color{red}$x+y$}
\psfrag{x-y}[l]{\color{red}$x-y$}
\includegraphics[height=55mm]{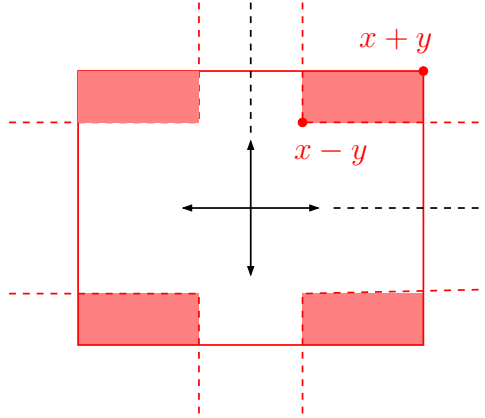}
\caption{Support of \,$\gamma_{x,y}$ \,for the root system \,$A_1\!\times\!A_1$}
\label{fig:supportA1A1}
\end{center}
\end{figure}

\begin{figure}[ht]
\begin{center}
\psfrag{x+y}[c]{\color{red}$x+y$}
\psfrag{x-y}[c]{\color{red}$x-y$}
\hspace{8mm}\includegraphics[height=75mm]{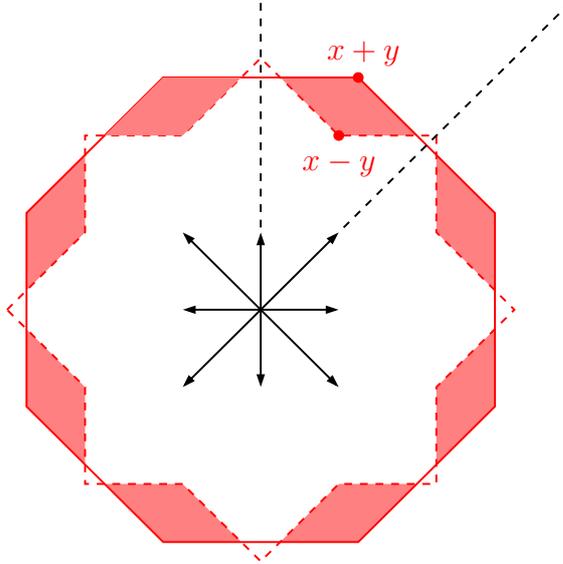}
\caption{Support of \,$\gamma_{x,y}$ \,for the root system \,$B_2$}
\label{fig:supportB2}
\end{center}
\end{figure}

\begin{proof}
Let $h\!\in\!\mathcal{C}_c^\infty(\mathbb{R}^N)$
be an auxiliary radial function such that
$$
\int_{\mathbb{R}^N}\hspace{-1mm}h(x)\,w(x)\,dx=1
$$
and \,$\supp h\!\subset\!-\co(G.u)$,
where $u\!\in\!\Gamma_{\!+}$ is a unit vector.
For every \,$\varepsilon\!>0$ \,and \,$x,y,z\!\in\!\mathbb{R}^N$,
\,set
$$
\gamma_\varepsilon(x,y,z)\,
=\,{\textstyle\frac1{c^2}}\int_{\mathbb{R}^N}\hspace{-1mm}
\mathcal{D}h(\varepsilon\,\xi)\,
E(i\,\xi,x)\,E(i\,\xi,y)\,E_k(-i\,\xi,z)\,w(\xi)\,d\xi\,.
$$
Firstly, according to \eqref{PaleyWiener0} and \eqref{EV},
$$
\xi\,\longmapsto\,\mathcal{D}h(\varepsilon\,\xi)\,E(i\,\xi,x)\,E(i\,\xi,y)
$$
is an entire function on \,$\mathbb{C}^N$ satisfying
\begin{equation}\label{PW}
\textstyle
\bigl|\,\mathcal{D}h(\varepsilon\,\xi)\,E(i\,\xi,x)\,E(i\,\xi,y)\,\bigr|\;
\le\,C_M\,(1\!+\!|\xi|)^{-M}\,
e^{-\langle\,g_0.(x_+\hspace{-.25mm}+\hspace{.25mm}y_+\hspace{-.25mm}
+\hspace{.25mm}\varepsilon\hspace{.2mm}u),\,(\Im\xi)_+\rangle}\,,
\end{equation}
where $g_0$ is the longest element in $G$,
which interchanges the chambers \,$\Gamma_{\!+}$ and $-\Gamma_{\!+}$\,.
Secondly,
\vspace{-2mm}
\begin{equation*}
\begin{aligned}
\langle\hspace{.5mm}\gamma_{x,y},f\hspace{.5mm}\rangle\hspace{.5mm}
&=\,{\textstyle\frac1c}
\int_{\mathbb{R}^N}\hspace{-1mm}\mathcal{D}f(\xi)\,
E(i\,\xi,x)\,E(i\,\xi,y)\,w(\xi)\,d\xi\\
&=\,{\textstyle\lim_{\,\varepsilon\to0}\,\frac1c}
\int_{\mathbb{R}^N}\hspace{-1mm}
\mathcal{D}h(\varepsilon\,\xi)\,\mathcal{D}f(\xi)\,
E(i\,\xi,x)\,E(i\,\xi,y)\,w(\xi)\,d\xi\\
&=\,{\textstyle\lim_{\,\varepsilon\to0}}
\int_{\mathbb{R}^N}\hspace{-1mm}f(z)\,\gamma_\varepsilon(x,y,z)\,w(z)\,dz
\end{aligned}
\end{equation*}
i.e. the distribution \,$\gamma_{x,y}$
\hspace{.25mm}is the weak limit of the measures
\,$\gamma_\varepsilon(x,y,z)\,w(z)\,dz$\,.
Thirdly, notice the symmetries
\begin{equation}\label{symmetry}
\gamma_\varepsilon(x,y,z)\,=\,\begin{cases}
\,\gamma_\varepsilon(y,x,z)\,,\\
\,\gamma_\varepsilon(g.x,g.y,g.z)
\quad\forall\;g\!\in\!G\cup\{-\text{Id}\}\,,\\
\,\gamma_\varepsilon(-z,y,-x)=\gamma_\varepsilon(x,-z,-y)\,.
\end{cases}
\end{equation}

If $G$ is crystallogaphic, we use Corollary \ref{PaleyWiener2}
(actually the third version of the Paley--Wiener theorem in \cite{J2}),
and deduce from \eqref{PW} that the function 
\,$z\,\longmapsto\gamma_\varepsilon(x,y,z)$
\,is supported in
\begin{equation*}
\co\,\{G.(x_+\hspace{-.75mm}+\hspace{-.25mm}y_+
\hspace{-.75mm}+\hspace{-.25mm}\varepsilon\,u)\}
=\hspace{.5mm}\co(G.x)+\co(G.y)+\varepsilon\co(G.u)\,.
\end{equation*}
Equivalently,
\begin{equation*}
\gamma_\varepsilon(x,y,z)\ne0
\quad\Longrightarrow\quad
z_+\!\prec
x_+\hspace{-.75mm}+\hspace{-.25mm}y_+
\hspace{-.75mm}+\hspace{-.25mm}\varepsilon\hspace{.25mm}u\,.
\end{equation*}
\vspace{-3mm}

\noindent
Using the symmetries \eqref{symmetry},
we see that \,$\gamma_\varepsilon(x,y,z)\ne0$ \,implies also
\begin{equation*}
\begin{cases}
\,-\hspace{.25mm}g_0.x_+\!\prec
-\hspace{.25mm}g_0.z_+\hspace{-.75mm}+\hspace{-.25mm}y_+
\hspace{-.75mm}+\hspace{-.25mm}\varepsilon u
&\text{i.e.}\quad
z_+\!\succ
x_+\hspace{-.75mm}+\hspace{-.25mm}g_0.y_+
\hspace{-.75mm}+\hspace{-.25mm}\varepsilon\hspace{.25mm}g_0.u\,,\\
\,-\hspace{.25mm}g_0.y_+\!\prec
-\hspace{.25mm}g_0.z_+\hspace{-.75mm}+\hspace{-.25mm}x_+
\hspace{-.75mm}+\hspace{-.25mm}\varepsilon u
&\text{i.e.}\quad
z_+\!\succ
g_0.x_+\hspace{-.75mm}+\hspace{-.25mm}y_+
\hspace{-.75mm}+\hspace{-.25mm}\varepsilon\hspace{.25mm}g_0.u\,.
\end{cases}
\end{equation*}
The conclusion of Theorem \ref{support} in the crystallographic case
is obtained by letting \,$\varepsilon\!\to\!0$\,.

If $G$ it not crystallographic,
we can only use the spherical Paley--Wiener theorem
and we obtain this way that
\,$\gamma_\varepsilon(x,y,z)\ne0$ \,implies
\begin{equation*}
\begin{cases}
\,|z|\le|x|+|y|+\varepsilon\,,\\
\,|x|\le|z|+|y|+\varepsilon\,,\\
\,|y|\le|x|+|z|+\varepsilon\,,
\end{cases}
\end{equation*}
hence
\begin{equation*}
\bigl|\,|x|-|y|\,\bigr|-\varepsilon\le|z|\le|x|+|y|+\varepsilon\,.
\end{equation*}
We conclude again by letting \,$\varepsilon\!\to\!0$\,.
\end{proof}

\end{document}